\documentclass{article}

\usepackage{amsmath,amsthm,amssymb}
\usepackage{amsfonts}
\usepackage{rotating}
\usepackage{euscript}
\usepackage{pst-node}
\usepackage{epsfig}

\def\be{\begin{equation}}
\def\ee{\end{equation}}

\def\C{{\mathbb C}} 
\def\N{{\mathbb N}} 
\def\P{{\mathbb P}}

\def\R{{\mathbb R}}

\def\phi{{\varphi}}

\def\cos{{\rm cos\,}}

\def\mod{{\rm mod\ }}
\def\tilde{\widetilde}

\def\bp{\begin{proposition}}
\def\ep{\end{proposition}}

\def\bt{\begin{theorem}}
\def\et{\end{theorem}}
\def\br{\begin{remark}}
\def\er{\end{remark}}
\def\be{\begin{equation}}
\def\bee{\begin{equation*}}
\def\la{\label}

\def\ee{\end{equation}}
\def\eee{\end{equation*}}
\def\bl{\begin{lemma}}
\def\el{\end{lemma}}
\def\bc{\begin{corollary}}
\def\ec{\end{corollary}}
\def\pr{\noindent{\it Proof. }}

\def\bd{\begin{definition}}
\def\ed{\end{definition}}
\input epsf.sty

\newtheorem{theorem}{Theorem}[section]
\newtheorem{lemma}{Lemma}[section]
\newtheorem{definition}{Definition}[section]
\newtheorem{corollary}{Corollary}[section]
\newtheorem{proposition}{Proposition}[section]
\newtheorem{remark}{Remark}[section]

\begin{document}
\title{Weak and strong composition conditions for the Abel differential equation}
\author{F. Pakovich}
\date{}

\maketitle

\begin{abstract}  
We establish an equivalence between two forms of the composition condition for the Abel differential equation with trigonometric coefficients.
\end{abstract}

\section{Introduction} Let  $\R_t[\theta]$ be the ring of trigonometric polynomials over $\R$, that is 
the ring generated over $\R$ by the functions 
$\cos \theta$, $\sin \theta$.
The center problem for the Abel differential equation 
\be \label{a} \frac{dr}{d\theta}=\widehat l(\theta)r^3+\widehat m(\theta)r^2,\ee
where $\widehat l,$ $\widehat m\in \R_t[\theta],$
is to find conditions implying that 
all its solutions are periodic on $[0,2\pi]$ whenever the initial condition is small enough. 
This problem is of a great interest because of its relation with 
the  classical Poincar\'e center-focus problem about the characterization 
of planar vector fields
\be \la{sys1}
\left\{
\begin{array}{rcl}
\dot{x} & =& -y+F(x,y), \\ 
\dot{y} & =&x+G(x,y),\\
\end{array} 
\right.
\ee
where $F(x,y),$ $G(x,y)$ are polynomials without constant and linear terms,
whose integral trajectories 
are closed in a neighborhood of the origin.
Namely, it was shown in \cite{cher} that in the case where $F(x,y),$ $G(x,y)$ are homogeneous and of the same degree, the Poincare problem reduces to the center problem for 
Abel equation \eqref{a}. The center problem for the Abel equation and its modifications are the subject 
of many recent papers involving different approaches and techniques (see e. g. \cite{abc}, \cite{alw1}, \cite{bpy}, \cite{bry}, \cite{bryd},  \cite{ga1}, \cite{i}, \cite{f}, \cite{ggl} and the bibliography therein).

Set
\be \label{form} l(\theta)=\int_0^{\theta}\widehat l(s)ds , \ \ \ m(\theta)=\int_0^{\theta}\widehat m(s)ds.\ee 
The following ``composition condition'' introduced in \cite{al} is  sufficient 
for equation \eqref{a} to have a center:
there exist $C^1$-functions $\tilde l,\tilde m, w$ with $w$ being $2\pi$-periodic such that  
\be \label{c} l(\theta)=\tilde l(w(\theta)), \ \ \ \ m(\theta)=\tilde m(w(\theta)). \ee
Indeed, if \eqref{c} holds, then any solution of \eqref{a} has the form $y(\theta)=\tilde y(w(\theta)),$ where 
$\tilde y$ is a solution of the equation 
$$ \frac{dr}{d\theta}=\tilde l^{\prime}(\theta)r^3+\tilde m^{\prime}(\theta)r^2,$$
implying that $y(0)=y(2\pi).$

In general, the composition condition is not necessary for \eqref{a} to have a center (\cite{alw}, \cite{alw1}). However, the composition condition 
is necessary and sufficient for some stronger forms of the center condition as well as for some other conditions related to the center problem
(see e. g. \cite{ga1}, \cite{ggl}). 
In fact, in all such cases the following apparently stronger condition imposed on $l$ and $m$ is satisfied:  
there exist a {\it trigonometric polynomial} $w$  and {\it polynomials} $\tilde l, \tilde m$ such that equalities 
\eqref{c} hold.
In this note we show that the last conditions is actually equivalent to the composition condition. More precisely, we prove the following statement:

\bt \label{mai} Let $l,m\in \R_t[\theta].$ Assume that there exist continuous functions $\tilde l,\tilde m, w$ with $w$ being $2\pi$-periodic 
such that 
the equalities  
$$ l(\theta)=\tilde l(w(\theta)), \ \ \ \ m(\theta)=\tilde m(w(\theta))$$ hold. 
Then they hold for some $\tilde l,$ $\tilde m\in \R[x]$ and $w\in \R_t[\theta]$.
\et

Thus, despite its analytic nature the composition condition turns out to be essentially algebraic. In particular, it can be expressed in terms of algebraic conditions imposed on coefficients of corresponding trigonometric polynomials.

\section{Proof of Theorem 1.1} Denote by $\R_t(\theta)$ the quotient field of $\R_t[\theta]$. It is well known that $\R_t(\theta)$ is isomorphic to 
the field $\R(x)$, where the isomorphism $\psi:\,\R_t(\theta)\rightarrow\R(x)$ is given by formulas 
$$\psi(\sin \theta)=\frac{2x}{1+x^2}, \ \ \ \psi(\cos \theta)=\frac{1-x^2}{1+x^2}, \ \ \ \psi^{-1}(x)=\tan\left(\frac{\theta}{2}\right).$$ In particular, this implies by the 
L\"uroth theorem that any subfield $k$ of $\R_t(\theta)$ has the form $k=\R(b)$ for some $b\in \R_t(\theta).$

\bl \label{l1} Let $l,m$ be non-constant trigonometric polynomials.  Assume that there exist continuous functions $\tilde l, \tilde m, w$ such that equalities \eqref{c} hold.
Then the field $\R(l,m)$ is distinct from the field $\R(\tan\left(\frac{n\theta}{2}\right))$ for any $n\geq 1$.
\el 
\pr Assume that $\R(l,m)=\R(\tan\left(\frac{n\theta}{2}\right))$ for some $n\geq 1.$
Then
there exist $u\in \R(x,y)$ such that 
\be \label{ccc} \tan\left(\frac{n\theta}{2}\right)=u(l,m).\ee
Clearly, conditions \eqref{c} and \eqref{ccc} imply that for any  $\theta_1, \theta_2\in \R$ the equality 
\be \label{tg} \tan\left(\frac{n\theta_1}{2}\right)=\tan\left(\frac{n\theta_2}{2}\right)\ee holds whenever \be \label{zxc} w(\theta_1)=w(\theta_2).\ee
On the other hand, equality \eqref{tg} holds if and only if $$\theta_1-\theta_2\equiv 0 \ \mod \frac{2\pi}{n}.$$ 
Therefore, in order to prove the lemma it is enough to find $\theta_1, \theta_2\in \R$ such that \eqref{zxc} holds but 
\be \label{tg2}\theta_1-\theta_2\not\equiv 0 \ \mod \frac{2\pi}{n}.
\ee

Since the function $w$ is continuous and $2\pi$-periodic, it attains its maximum value $x_0$ on $\R.$ Furthermore, 
it follows easily from the intermediate value theorem that 
for any positive $\epsilon$ which is small enough the equation 
$w(\theta)=x_0-\epsilon$ has at least two distinct roots $\theta_{1},$ $\theta_{2}$ 
which satisfy \eqref{tg2}.
\qed

\vskip 0.2cm

The following lemma, describing subfields of $\R_t(\theta)$ containing trigonometric polynomials,
is proved in the paper \cite{ggl} (Proposition 21)  and in the paper \cite{ga1} (Theorem 5). However, the proofs given in \cite{ggl}, \cite{ga1} are quite complicated and occupy several pages. Below we provide a short independent proof which is based on the fact that the ring $\R_t[\theta]$ is isomorphic to a subring of the ring $\C[z,1/z]$ of complex Laurent polynomials, where an isomorphism $\phi:\,\R_t[\theta]\,\rightarrow \C[z,1/z]$ is given  by the formulas: 
\be \la{t} \cos \theta\rightarrow \left(\frac{z+1/z}{2}\right), \ \ \ 
\sin \theta\rightarrow \left(\frac{z-1/z}{2i}\right). \ee 
Notice that the isomorphism $\phi$ can be used for a construction of a comprehensive decomposition theory of trigonometric polynomials (see \cite{ppaa}).

\vskip 0.2cm

\bl \label{l2} Let $k$ be a subfield of $\R_t(\theta)$ containing a non-constant trigonometric polynomial. Then either $k=\R(\tan(\frac{n\theta}{2}))$ for some $n\in \N,$ or $k=\R(b)$ for some trigonometric polynomial $b.$
\el 
\pr
For brevity,  we will denote the ring $\C[z,1/z]$ by ${\cal L}[z]$ and  the image of 
$\R_t[\theta]$ in  $\cal L$ under the isomorphism $\phi$ by ${\cal L_{\R}}[z]$. It is easy to see that ${\cal L_{\R}}[z]$ consists of Laurent polynomials 
$L$ such that $\bar L(1/z)=L(z),$ where $\bar L$ denotes the Laurent polynomial obtained from $L$ by complex conjugation of all its coefficients. The isomorphism $\phi$ extends to an isomorphism between the quotient field $\R_t(\theta)$  of $\R_t[\theta]$
and the  quotient field  ${\cal L_{\R}}(z)$ of ${\cal L_{\R}}[z]$. Clearly, the field ${\cal L_{\R}}(z)$
consists of rational functions $R$ satisfying the equality 
$\bar R(1/z)=R(z).$

Assume that  $k$ is a subfield of $\R_t(\theta)$ containing a non-constant trigonometric polynomial $l$. Let $b$ be an element of $\R_t(\theta)$
such that $k=\R(b)$ and $A\in \R(x)$ be a rational function such that $l(\theta)=A(b(\theta)).$  
Set $L=\phi(l),$ $B=\phi(b)$. Clearly, 
$L(z)=A(B(z))$. Further, since $L$ is a Laurent polynomial we have:
$$L^{-1}\{\infty \}=B^{-1}\{A^{-1}\{\infty \}\}= \{0,\infty\},$$ implying that
the set $A^{-1}\{\infty \}$ contains at most two points. In more details, either 
$$A^{-1}\{\infty \}=\{a\} \ \ {\rm and} \ \ B^{-1}\{a \}=\{0,\infty\},$$ for some $a\in \C\P^1,$ or 
$$A^{-1}\{\infty \}=\{a,b\} \ \  {\rm and} \ \  B^{-1}\{a,b \}=\{0,\infty\},$$ for some $a, b\in \C\P^1.$

It is easy to see that in the first case 
there exists a rational function 
$\mu\in \C(z)$ of degree one such that $A(\mu (z))$ is a polynomial and $\mu^{-1}(B(z))$ is a Laurent polynomial, while in the second case there exists
a rational function $\mu\in \C(z)$ of degree one such that $A(\mu (z))$ is a Laurent polynomial and $\mu^{-1}(B(z))=z^d,$ $d>1$.

Since any polynomial with real coefficients is a product of linear and quadratic polynomials with real coefficients, in the first case
the equality $A^{-1}\{\infty \}=\{a\}$ 
implies that $a\in \R$, unless $a=\infty.$ Therefore, setting $\mu(z)=a+1/z$, 
we can assume that $\mu$ has {\it real} coefficients. Since $B\in {\cal L_{\R}}(z)$, this implies that 
the function $\mu^{-1}(B(z))$ is contained ${\cal L_{\R}}[z]$, and hence  
$\mu^{-1}(b(\theta))$ is a 
trigonometric polynomial, since 
$\phi$ is an isomorphism. Clearly, 
this polynomial generates the field $k.$

In the second case,  composing $\mu$ 
with an other rational function of degree one, we obtain a rational function $\mu_1\in \C(z)$ of degree one such that 
$$\mu^{-1}_1(B(z))=\frac{1}{i}\frac{z^d-1}{z^d+1}=\frac{1}{i}\left(\frac{z^{d/2}-z^{-d/2}}{z^{d/2}+z^{-d/2}}\right)=\phi(\tan(d\theta/2)).$$
Since the rational functions $\phi(\tan(d\theta/2))$ and $B(z)$ are contained in ${\cal L_{\R}}(z)$, the last equality implies  easily that 
$\bar \mu^{-1}_1=\mu^{-1}_1.$
Therefore, $\mu^{-1}_1\in \R(x)$ and 
$\mu^{-1}_1(b)=\tan(d\theta/2).$ \qed

\vskip 0.2cm

Theorem 1.1 follows from the above lemmas. Indeed, 
by Lemma \ref{l1}, the field $k=\R(l,m)$  is distinct from the field $\R(\tan\left(\frac{n\theta}{2}\right))$ for any $n\geq 1$. Therefore, 
by Lemma \ref{l2} this field is generated by some trigonometric polynomial $w$
implying that equalities \eqref{c} hold for some $\tilde l,$ $\tilde m\in \R[x]$ and $w\in \R_t[\theta]$.

\end{document}